\documentclass[12pt, leqno, oneside]{amsart}
\usepackage{enumerate}
\usepackage{enumerate}
\newtheorem{theorem}{Theorem}

\newtheorem{lemma}[theorem]{Lemma}

\newtheorem{example}[theorem]{Example}
\theoremstyle{definition}
\newtheorem{definition}[theorem]{Definition}
\theoremstyle{remark}
\newtheorem{remark}[theorem]{Remark}

\def\aut{\operatorname{Aut}}

\def\tr{\operatorname{Tr}}
\def\prop{\operatorname{Prop}}
\title[Serre problem]
{Serre problem for unbounded pseudoconvex Reinhardt domains in $\mathbb C^2$}
\author[\L.~Kosi\'nski]
{\L ukasz Kosi\'nski}
\address{Instytut Matematyki\\ Uniwersytet Jagiello\'nski\\ \L ojasiewicza 6,
30-348 Krak\'ow\\ Poland} \email{lukasz.kosinski@gazeta.pl}  


\thanks{Research partially supported by the KBN grant N$^{\text o}$ \textbf{N N201 271435} and by the foundation of A.~Krzy\.zanowski}

\begin{document}
\maketitle
\begin{abstract}
We give a characterization of non-hyperbolic pseudoconvex Reinhardt domains in $\mathbb C^2$ for which the answer to the Serre problem is positive. Moreover, all non-hyperbolic pseudoconvex Reinhardt domains in $\mathbb C^2$ with non-compact automorphism group are explicitly described.
\end{abstract}  
\section{Statement of result}
Throughout this paper the class of Stein domains $D$ for which the answer to the Serre problem is positive (with the fiber equal to $D$) is denoted by $\mathfrak S$, i.e. $D\in\mathfrak S$ if for any Stein manifold $B,$ any holomorphic fiber bundle $E\to B$ with the base $B$ and the Stein fiber $D$ is Stein. In 1953 J-P. Serre raised a question whether all Stein manifolds are in $\mathfrak S.$

Despite many positive results the answer to the Serre problem is in general negative. Actually, Skoda proved that $\mathbb C^2\notin\mathfrak S$ (see \cite{Sko}). The counterexamples with bounded domains as fibers were found by G. Coeur\'e and J.-J.~Loeb (\cite{Coe}).
In \cite{pzwo} P.~Pflug and W.~Zwonek gave a characterization of all hyperbolic Reinhardt domains of $\mathbb C^2$ not in $\mathfrak S.$ 
Next K.~Oeljeklaus and D. Zaffran solved the Serre problem for bounded Reinhardt domains in $\mathbb C^3.$ 
Recently a classification result for bounded Reinhardt domains of $\mathbb C_*^d$ for arbitrary $d\geq 2$ has been obtained by D. Zaffran in \cite{Zaf}. In particular, the case $\mathbb C_*^d,$ $d\geq 2,$ may be deduced from this result (it follows easily from \cite{Zaf}, Main Theorem).

In the paper we deal with non-hyperbolic Reinhardt domains in $\mathbb C^2$ and we solve the Serre problem for them.
The main goal is to show the following
\begin{theorem}\label{main} Let $D$ be a pseudoconvex non-hyperbolic Reinhardt domain. Then $D\notin\mathfrak S$ if and only if 
$\mathbb C_*^2\subset D$ or $D$ is algebraically equivalent to a domain of the form 
\begin{equation} \{(z_1,z_2)\in\mathbb C_*^2:\ |z_1||z_2|^{p\pm\sqrt{q}}<1\},\end{equation} where $p,q\in\mathbb Q,\ q>0,\ \sqrt{q}\in\mathbb R\setminus\mathbb Q.$
\end{theorem}

\begin{remark}\label{uw}
Note that for the pseudoconvex Reinhardt domain $D$  the condition $\mathbb C_*^2\subset D$ means that $D\in\{\mathbb C^2,\mathbb C_*^2,\mathbb C\times\mathbb C_*,\mathbb C_*\times\mathbb C\}.$
\end{remark}

The following result gives a description of non-hyperbolic pseudoconvex Reinhardt domains in $\mathbb C^2$ whose group of automorphisms is non-compact. 

\begin{theorem}\label{compact}Let $D$ be a non-hyperbolic pseudoconvex Reinhardt domain in $\mathbb C^2.$ Then the group $\aut(D)$ is non-compact
if and only if the logarithmic image of the domain $D$ contains an affine line or (up to a permutation on components) $D$ is contained in $\mathbb C\times\mathbb C_*,$ $\{0\}\times\mathbb C_*\subset D$ and $\log D=\{(t,s):\ t<\psi(s),\ s\in\mathbb R\},$ where $\psi:\mathbb R\to \mathbb R$ is a concave function satisfying the property $\psi(\beta+s)-\psi(s)=\alpha+ks,\ s\in\mathbb R,$ for some $\alpha,\beta\in\mathbb R,\ \beta\neq 0,\ k\in\mathbb Z_*.$
\end{theorem}

\begin{remark}
The examples of a functions $\psi$ appearing in Theorem~\ref{compact} are the functions $\psi(s)=as^2+bs+c,$ $s\in\mathbb R,$ where $a<0,$ $b,c\in\mathbb R.$\end{remark}

\bigskip

The author expresses his thanks to the referee for his valuable suggestions to improve the shape of the original paper. Also, I would like to express my gratitude to professor W\l odzimierz Zwonek for helpful discussions.

\section{Preliminaries}
Here is some notation. By $\mathbb D$ we denote the open unit disc in the complex plane. Let $\mathbb A(r,R)=\{z\in\mathbb C:\ r<|z|<R\},$ $-\infty<r<R<\infty,$ $R>0.$ Note that if $r>0,$ then $\mathbb A(r,R)$ is an annulus. For simplicity we put $\mathbb A(r)=\mathbb A(1/r,r),$ $r>0.$ Moreover, for a domain $D$ in $\mathbb C^n,$ the set $D\setminus\{0\}$ is denoted by $D_*.$  

For a Reinhardt domain $D\subset\mathbb C^n$ we define $$\log D=\{(x_1,\ldots,x_n)\in\mathbb R^n:\ (e^{x_1},\ldots,e^{x_n})\in D\}.$$ It is well known that, for any pseudoconvex domain in $\mathbb C^n,$ its logarithmic image is convex.

Put $$\mathbb C(\lambda):=\begin{cases}\mathbb C,\ \text{if}\quad \lambda\geq0,\\ \mathbb C_*,\ \text{if}\quad \lambda<0,\end{cases}$$ and $$\mathbb C(\alpha):=\mathbb C(\alpha_1)\times\ldots\times \mathbb C(\alpha_n),\quad \alpha\in (\mathbb R^n).$$

Given $\alpha\in\mathbb R^n$ and $z\in\mathbb C(\alpha),$ we put $z^{\alpha}=z_1^{\alpha_1}\ldots z_n^{\alpha_n}$ \footnote{We define $0^0:=1$.}. Moreover, we define $$\varphi_A(z):=z^A:=(z^{A^1},\ldots,z^{A^m}),\quad z\in\mathbb C(A):=\bigcap_{j=1}^m\mathbb C(A^j),$$ where $A= (A_k^j)_{j=1,\ldots,m,k=1,\ldots,n} \in\mathbb Z^{m\times n}.$

\medskip

Non-hyperbolic Reinhardt domains appearing in Theorem~\ref{main} have been characterized in the following way:
\begin{theorem}[see~\cite{Zwo}]\label{hz}
Let D be a pseudoconvex Reinhardt domain in $\mathbb C^n.$ Then the following conditions are equivalent:
\begin{itemize}
\item[-] D is (Kobayashi, Carath\'eodory or Brody) hyperbolic;
\item[-] D is algebraically equivalent to a bounded domain - in other words, there is an $A\in\mathbb Z^{n\times n}$ with $|\det A|=1$ such that $\varphi_A(D)$ is bounded and $(\varphi_A)|_D$ is a biholomorphism onto the image;
\item[-] $\log D$ contains no straight lines and $D\cap (\mathbb C^{j-1}\times\{0\}\times\mathbb C^{n-j})$ is either empty or hyperbolic (viewed as domains in $\mathbb C^{n-1}$), $j=1,\ldots,n.$
\end{itemize}
\end{theorem}

In view of Theorem~\ref{hz}, non-hyperbolic pseudoconvex Reinhardt domains in $\mathbb C^2$ not containing $\mathbb C_*^2$ may be divided into two classes. 

The first one contains non-hyperbolic, pseudoconvex Reinhardt domains $D$ such that $D\cap\mathbb C_*^2$ is hyperbolic.

The second one consists of domains whose logarithmic images contain a straight line. One may prove that such domains are algebraically equivalent to \begin{equation}\label{d}D_{\alpha,r_-,r_+}:=\{(z_1,z_2)\in\mathbb C(\alpha):\ r_-<|z_1|^{\alpha_1}|z_2|^{\alpha_2}<r_+\},\footnote{Note that replacing $(\alpha_1,\alpha_2)$ by $(-\alpha_1,\alpha_2)$ allows us to add or to remove the axis $\{0\}\times \mathbb C$ in the case when $\alpha_1\neq 0$. In other words $\{(z_1,z_2)\in\mathbb C_*\times \mathbb C(\alpha_2) :\ r_-<|z_1|^{\alpha_1}|z_2|^{\alpha_2}<r_+\},$ where $\alpha_1>0,$ $\alpha_2\in\mathbb R$ is also covered by this symbol (it is algebraically equivalent to $D_{(-\alpha_1,\alpha_2),r_-,r_+}$). If $\alpha_1=0,$ then $D_{\alpha,r_-,r_+}=\mathbb A(r_-,r_+)\times\mathbb C$. Moreover, $\mathbb A(r_-,r_+)\times\mathbb C_*$ is algebraically equivalent to $D_{(1,-1),r_-,r_+}$.}\end{equation} where $\alpha=(\alpha_1,\alpha_2)\in\mathbb R^2_*,$ $-\infty<r_-<r_+<\infty,$ $r_+>0$ .

\begin{definition}A domain $D_{\alpha,r_-,r_+}$ is called to be of \emph{rational type} if $t\alpha\in(\mathbb Q^2)_*$ for some $t>0.$
 
Otherwise $D_{\alpha,r_-,r_+}\subset \mathbb C^2$ is said to be of \emph{irrational type}.
\end{definition}

It is clear that $D_{\alpha,r^-,r^+}$ is algebraically equivalent to a domain of one of the following types: 
\begin{enumerate}
\item[(a)] $D_{\beta}:=\{(z_1,z_2)\in\mathbb C\times \mathbb C(\beta):\ |z_1||z_2|^{\beta}<1\},\ \beta\in\mathbb R,$
\item[(b)] $D_{\beta}^*:=\{(z_1,z_2)\in\mathbb C_*\times\mathbb C(\beta):\ 0<|z_1||z_2|^{\beta}<1\},\ \beta\in \mathbb R_{\geq 0},$
\item[(c)] $D_{\beta,r}:=\{(z_1,z_2)\in\mathbb C\times \mathbb C(\beta):\ 1/r<|z_1||z_2|^{\beta}<r\},\ \beta\in\mathbb R_{\geq 0},\ r>1.$
\end{enumerate}

\begin{remark}\label{rationalaut} Observe that any Reinhardt domain $D_{\alpha,r_-,r_+}$ of rational type is algebraically equivalent to one of the following domains:

\begin{enumerate} \item[(d)] $\mathcal D\times\mathcal C,$ where $\mathcal D$ is either a disc, a pointed disc or an annulus and $\mathcal C\in\{\mathbb C,\mathbb C_*\},$
\item[(e)] $\{(z,w)\in\mathbb C^2:\ |z^p||w^q|<1\},$ where $p,q\in\mathbb N$ are relatively prime,
\item[(f)] $\{(z,w)\in\mathbb C_*\times\mathbb C:\ |z^p||w^q|<1\},$ where $p,q$ are relatively prime natural numbers, $(p,q)\neq (1,1)$.
\end{enumerate}
\end{remark}
\begin{proof} We will show that any domain of the form  (a), (b), (c) with rational $\beta$ is algebraically equivalent to one of the enumerated above domains. 

Domains of the form $D_0,$ $D_0^*$ and $D_{0,r}$ are in (d). 

Note that domains of the form (e) may be written as $D_{\frac{q}{p}}$. Moreover, domains of the form (f) are algebraically equivalent to $D_{-\frac{q}{p}}.$ Therefore any domain of the form (a) with $\beta\neq 0$ is of the form (e) or (f).

Fix a rational $\beta>0$ and assume that $\beta=p/q,$ where $p$ and $q$ are relatively prime natural numbers. There are $m,n\in\mathbb Z$ such that $pm+qn=1.$ Observe that the mapping $D_{\beta}^*\ni (z_1,z_2)\to (z_1^qz_2^p,z_1^{-m}z_2^n)\in \mathbb D_*\times \mathbb C_*$ is biholomorphism.

Similarly, one may show that $D_{\beta,r}$ is algebraically equivalent to $\mathbb A(r)\times\mathbb C_*.$

What remains to do is to observe that the set $\{(z,w)\in \mathbb C_*\times \mathbb C:\ |z||w|<1\}$ is algebraically equivalent to $\mathbb D\times \mathbb C_*.$
\end{proof}

\bigskip
Our paper is organized as follows. In Section~\ref{sek1} we give the proof of Theorem~\ref{main} for Reinhardt domains of the form (\ref{d}) of the rational type. 

Section~4 is devoted to solving the Serre problem for Reinhardt domains of the form (\ref{d}) of irrational type. We obtain there a  natural correspondence between the automorphisms of Reinhardt domains of irrational type and the famous \textit{Pell's equation}. 

In Section~\ref{sek3} we proof Theorem~\ref{main} for non-hyperbolic Reinhardt domains which are hyperbolic after removing the axes $(\mathbb C\times\{0\})\cup(\{0\}\times\mathbb C).$

It seems to us probable that the solution of the Serre problem for $\mathbb C\times\mathbb C_*$ is known, however we could not find it in the literature, therefore in the last section we present ideas explaining how to easily extend the procedure used in \cite{Dem} to this domain. 

Finally, in Section~\ref{pomt} using obtained results we present the proofs of Theorem~\ref{main} and Theorem~\ref{compact}.

\bigskip

Two following theorems providing us with some classes of domains for which the answer to the Serre problem is positive will be useful. The first of them is the so-called \textit{Stehl\'e criterion}.
\begin{theorem}[see \cite{Mok,Stel}]\label{Ste} Let $D$ be a domain in $\mathbb C^n.$ If there exists a real-valued plurisubharmonic exhaustion function $u$ on $D$ such that $u\circ f-u$ is bounded from above for any $f\in\aut(D),$ then $D\in\mathfrak S.$
\end{theorem}
Observe the Stehl\'e criterion is satisfied among others by domains whose group of automorphisms is compact.
\begin{theorem}[see \cite{Mok0}]\label{Mok0} Any open Riemann surface belongs to $\mathfrak S.$
\end{theorem}

Let us recall basic facts and definitions related to holomorphic fiber bundles which will be useful in the sequel. For more information we refer the reader to \cite{Her}.
 
Let $E$ be an arbitrary holomorphic fiber bundle with the fiber $X.$ The automorphisms of the form $\tau_{\alpha,\beta}=\tau_{\alpha}\circ \tau^{-1}_{\beta}\in\aut((\Omega_{\alpha}\cap \Omega_{\beta})\times X),$ where $\tau_{\alpha},$ $\tau_{\beta}$ are trivializations (with associated domains $\Omega_{\alpha}\times X, \Omega_{\beta}\times X,$ respectively) are so-called \textit{transition functions}. It is clear that \begin{equation}\label{com}\tau_{\alpha,\beta}\circ\tau_{\beta,\gamma}=\tau_{\alpha,\gamma}\quad \text{on}\quad (\Omega_{\alpha}\cap\Omega_{\beta}\cap\Omega_{\gamma})\times X.\end{equation}

On the other hand, having domains $(\Omega_{\alpha})_{\alpha}$ and a family of functions $(\tau_{\alpha,\beta})_{\alpha,\beta}\in\aut((\Omega_{\alpha}\cap\Omega_{\beta})\times X)$ satisfying the condition (\ref{com}) we may define a holomorphic fiber bundle with the base $B=\bigcup_{\alpha}\Omega_{\alpha}$ and the fiber $X$ by putting \begin{equation}\label{form}E=\left(\bigsqcup_{\alpha}(\Omega_{\alpha}\times X)\right)/\sim,\end{equation} gluing the charts $\Omega_{\alpha}\times X$ via identification: $(x_{\alpha},z_{\alpha})\sim (x_{\beta},z_{\beta})$ if and only if $x_{\alpha}=x_{\beta}\in \Omega_{\alpha}\cap \Omega_{\beta}$ and $\tau_{\alpha,\beta}(x_{\beta},z_{\beta})=(x_{\alpha},z_{\alpha}).$

Throughout this paper all holomorphic fiber bundles are understood to be of the form (\ref{form}).

\section{Rational case}\label{sek1}
We start this section with collecting the formulas for automorphisms of Reinhardt domains of rational type. Some of them may be found in~\cite{Jar-Pfl} and \cite{Sol}. They may be also derived from the description of proper holomorphic mappings between Reinhardt domains of rational type obtained in \cite{prop} and \cite{ja}. As mentioned in Remark~\ref{rationalaut}, any rational Reinhardt domain is algebraically equivalent to a domain of the form (d), (e) or (f).
\begin{theorem}\label{wymaut}$ $ \begin{enumerate}[(1)]
\item The group of automorphism of the domain $\{(z_1,z_2)\in \mathbb C^2:\ |z_1z_2|<1\}$ consists of the mappings $$(z_1,z_2)\to (z_1 f(z_1z_2), e^{i\theta} z_2/f(z_1z_2))$$ and $$(z_1,z_2)\to (z_2f(z_1 z_2), e^{i\theta} z_1/f(z_1z_2)),$$ where $f\in \mathcal O^*(\mathbb D)$ and $\theta \in \mathbb R.$

\item The group of automorphism of $\{(z_1,z_2)\in \mathbb C^2:\ |z_1^p||z_2^q|<1\},$ where $p$ and $q$ are relatively prime natural numbers, $(p,q)\neq (1,1)$, consists of the mappings of the form $$(z_1,z_2)\to \left(a^{1/p}(z_1^pz_2^q)z,e^{i\theta}\frac{z_2}{a^{1/q}(z_1^pz_2^q)}\right),$$ where $a\in \mathcal O^*(\mathbb D),$ $\theta\in \mathbb R.$

\item The group of automorphism of the domain $\{(z_1,z_2)\in \mathbb C_* \times \mathbb C:\ |z_1^p||z_2^q|<1\},$ where $p$ and $q$ are relatively prime natural numbers, $p\geq 1,$ $q\geq 2,$ consists of the mappings of the form $$(z_1,z_2)\to \left(a^{1/p} (z_1^pz_2^q)z_1, e^{i\theta} \frac{z_2}{a^{1/q} (z_1^pz_2^q)} \right),$$ where $a\in \mathcal O^*(\mathbb D),$ $\theta\in \mathbb R.$

\item Let $\mathcal D\subset \mathbb C$ be bounded and $\mathcal C\in \{\mathbb C,\mathbb C_*\}.$ The group of automorphisms of $\mathcal D\times \mathcal C$ consist of the mappings $$\mathcal D\times \mathcal C\ni (\lambda,z) \to (a(\lambda),f(\lambda,z)),$$ where $a\in \aut(\mathcal D)$. 

In the case when $\mathcal C=\mathbb C$ the mapping $f$ is of the form $$f(\lambda,z)=b(\lambda)z+c(\lambda),\quad \lambda\in \mathcal D,\ z\in \mathbb C,$$ where $b\in \mathcal O^*(\mathcal D)$ and $c\in \mathcal O(\mathcal D).$

If $\mathcal C=\mathbb C_*,$ then the mapping $f$ if of the form $$f(\lambda,z)=b(\lambda)z^{\epsilon},\quad \lambda \in \mathcal D,\ z\in \mathbb C_*,$$ where $b\in \mathcal O^*(\mathcal D)$ and $\epsilon=\pm1.$
\end{enumerate}
\end{theorem}

The following lemma due to Stehl\'{e} will be needed
\begin{lemma}[\cite{Stel}]\label{Stell} Let $B$ be a Stein manifold. Then for any covering $\{U_{\alpha}\}$ of $B$ there is a locally finite covering $\{B_j\}_{j\in\mathbb N}$ and a family of strictly plurisubharmonic, continuous functions $h_j$ in a neighborhood of $\overline{C_j},$ where $C_j=\bigcup_{i\leq j}B_i,$ such that
\begin{enumerate}
\item[(i)] each $B_j$ is relatively compact in some $U_{\alpha},$
\item[(ii)] $C_j$ is Stein for any $j,$
\item[(iii)] $h_j<0$ in $C_j\setminus B_j$ and $h_j>1$ in the neighborhood (in $C_j$) of $\overline{C_j\setminus C_{j-1}}\cap C_j,$ $j\in\mathbb N.$
\end{enumerate}
\end{lemma}

\begin{theorem}\label{dc}Any Reinhardt domain $D$ in $\mathbb C^2$ of the form (\ref{d}) of rational type belongs to the class $\mathfrak S.$
\end{theorem}

\begin{remark}Note that the function satisfying the Stehl\'e criterion  does not exist for any domain appearing in the previous theorem. As an example let us take $\mathbb D\times\mathbb C.$ It is sufficient to observe that for any $f\in\mathcal O(\mathbb D,\mathbb C_*)$ the mapping $$\mathbb D\times\mathbb C\ni(\lambda,z) \to (\lambda,f(\lambda)z) \in \mathbb D\times\mathbb C$$ is an automorphism of $\mathbb D\times\mathbb C,$ and the growth of the function $f$ may be arbitrarily fast.\end{remark}

\begin{proof}[Proof of Theorem~\ref{dc}] The proof is divided into four cases.

$(a)$ First we focus our attention on the case when the fiber is equal to $D=\{(z,w)\in\mathbb C^2:\ |z||w|<1\}.$

We will modify the idea used by Stehl\'e in the proof of Theorem~\ref{Ste}. Let $\pi:E\to B$ be a holomorphic fiber bundle with a Stein base $B$ and a fiber $D.$ As we have already mentioned, $E$ is assumed to be given by the formula (\ref{form}). 

Using the description of the group of automorphisms of the domain $D$ we get that for any $b\in\Omega_{\alpha}\cap\Omega_{\beta}$ a transition function $\tau_{\alpha,\beta}$ is of one of two following forms: $$\begin{cases}\tau_{\alpha,\beta}(b,z,w)=(b,zf_{\alpha,\beta}(b,zw),e^{i\theta_{\alpha,\beta}(b)}w/f_{\alpha,\beta}(b,zw))\quad \text{or}\\ \tau_{\alpha,\beta}(b,z,w)=(b,wf_{\alpha,\beta}(b,zw),e^{i\theta_{\alpha,\beta}(b)}z/f_{\alpha,\beta}(b,zw)),\end{cases}$$ ($(z,w)\in D$), for  some $f_{\alpha,\beta}(b,\cdot)\in\mathcal O^*(\mathbb D)$ and $\theta_{\alpha,\beta}(b)\in\mathbb R.$ Since $\tau_{\alpha,\beta}$ is holomorphic it follows that the mapping $\Omega_{\alpha}\cap\Omega_{\beta}\ni b\rightarrow e^{i\theta_{\alpha,\beta}(b)}\in\partial\mathbb D$ is holomorphic and therefore locally constant.

Let $\tilde\pi:\tilde E\to B$ denote a holomorphic fiber bundle with the base $B$ and the fiber $\mathbb D$ whose transition functions are defined by $\tilde\tau_{\alpha,\beta}(b,\lambda)=(b,e^{i\theta_{\alpha,\beta}(b)}\lambda),\ b\in\Omega_{\alpha}\cap\Omega_{\beta},$ $\lambda\in\mathbb D.$ It follows from Theorem~\ref{Mok0} that $\tilde E$ is Stein. Applying Lemma~\ref{Ste} to the family $\{\tilde\pi^{-1}(\Omega_{\alpha})\}$ we get a locally finite covering $\{B_j\}$ and a family of strictly plurisubharmonic continuous functions $\{h_j\}$ satisfying conditions (i)-(iii) of Lemma~\ref{Stell}. Let $$p:E\ni[(b,z,w)]\rightarrow[(b,zw)]\in\tilde E$$ be a natural surjection between $E$ and $\tilde E.$ 

For $j\in\mathbb N$ choose $\alpha_j$ such that $B_j$ is relatively compact in $\tilde\pi^{-1}(\Omega_{\alpha_j}).$ Let $\tau_j:=\tau_{\alpha_j}$ denote any trivialization of the fiber bundle $\pi:E\to B$ defined on $p^{-1}(\tilde\pi^{-1}(\Omega_{\alpha_j})).$ Put $$l(b,z,w)=\max\{ \log^+|z|,\log^+|w|\},\quad (b,z,w)\in B\times\mathbb C^2.$$ The choice of $\{B_j\}$ and a standard compactness argument guarantee that for any $i,j\in\mathbb N$
\begin{equation}\label{l}\sup\{l(\tau_i\circ\tau_j^{-1}(b,z,w))-l(b,z,w):\ [b,z,w]\in p^{-1}(B_i\cap B_j)\}<\infty.
\end{equation}

Put $$v_1(x)=\exp\Big(h_1\big(p(x)\big)+l\big(\tau_1(x)\big)\Big),\quad x\in p^{-1}(B_1).$$

Condition (\ref{l}) allows us to choose constants $d\in(0,1)$ and $M>0$ such that \begin{equation}dM\exp\Big(l\big(\tau_2(x)\big)\Big)\leq v_1(x)\leq 2M\exp\Big(l\big(\tau_2(x)\big)\Big),\quad  x\in p^{-1}(B_1\cap B_2).\end{equation} Define $$\tilde v_2(x)=2M\exp\Big(l\big(\tau_2(x)\big)\Big)\exp\Big(\big(1-h_2(p(x))\big)\log\frac{d}{2}\Big),\quad x\in p^{-1}(B_2).$$

It follows from the choice of $d,M$ that if $h_2(p(x))<0,$ $x\in p^{-1}(B_1\cap B_2),$ then $\tilde v_2(x)<v_1(x).$ Moreover, if $h_2(p(x))>1,$ $x\in p^{-1}(B_1 \cap B_2),$ then $\tilde v_2(x)>v_1(x).$ Therefore putting $$
v_2(x)=\begin{cases}v_1(x),\quad x\in p^{-1}(B_1\setminus B_2),\\ \max\{v_1(x),\tilde v_2(x)\},\quad x\in p^{-1}(B_1\cap B_2),\\ \tilde v_2(x),\quad x\in p^{-1}(B_2\setminus B_1),\end{cases}$$ we obtain a plurisubharmonic continuous function defined on $p^{-1}(B_1\cup B_2)=p^{-1}(C_2).$ 

Similarly, let $d'\in (0,1),$ $M'>0$ be such that \begin{equation}dM\exp\Big(l\big(\tau_3(x)\big)\Big)\leq v_2(x)\leq 2M\exp\Big(l\big(\tau_3(x)\big)\Big),\quad  x\in p^{-1}(B_3\cap C_2)\end{equation} Let us define $$\tilde v_3(x)=2M'\exp\Big(l\big(\tau_3(x)\big)\Big) \exp\Big(1-h_3\big(p(x)\big)\log\frac{d'}{2}\Big)\quad \text{for}\quad x\in p^{-1}(B_3),$$ and $$v_3(x)=\begin{cases}v_2(x),\quad x\in p^{-1}(C_2\setminus B_3),\\ \max\{v_2(x),\tilde v_3(x)\},\quad x\in p^{-1}(C_2\cap B_3),\\ \tilde v_3(x),\quad x\in p^{-1}(B_3\setminus C_2).\end{cases}$$ Exactly as before we check that $v_3\in\mathcal{PSH}(p^{-1}(C_3)).$ 

Repeating this procedure inductively one may obtain a sequence of functions $v_j\in\mathcal{PSH}(p^{-1}(C_j))$ such that $v_j\leq v_{j+1}$ on $p^{-1}(C_j)$ and $v_j=v_{j+1}$ on $p^{-1}(C_j\setminus B_{j+1}).$ 

Since the covering $\{B_j\}$ is locally finite, putting $v=\lim_j v_j$ we define properly a plurisubharmonic function on $E$.

Let $\tilde u$ be a strictly plurisubharmonic continuous exhaustion function on $\tilde E.$ It follows from the construction of the function $v$ that
\begin{equation}u=\max\{\tilde u \circ p,v \}
\end{equation} is a plurisubharmonic continuous exhaustion function on $E.$ 

Now, in order to show the Steinness of $E$ it suffices to repeat the argument used by N.~Mok in the proof of the improvement of the Stehl\'{e} criterion (see \cite{Mok}, Appendix III).

$(b)$ Let us consider the case when the fiber $D$ is of the form $D=\{(z,w)\in\mathbb C^2:\ |z|^p|w|^q<1\}$ for some natural relatively prime $p,q,\ (p,q)\neq(1,1).$ 

Let $\pi:E\to B$ be a holomorphic fiber bundle with the fiber $D.$ Using Theorem~\ref{wymaut} we infer that the transition functions are of the form
\begin{equation}\tau_{\alpha,\beta}(b,z,w)= \left(b,a_{\alpha,\beta}^{1/p}(b,z^pw^q)z,e^{i\theta_{\alpha,\beta}(b)} \frac{w} {a_{\alpha,\beta}^{1/q}(b,z^pw^q)} \right),\end{equation} where $a_{\alpha,\beta}\in\mathcal O^*((\Omega_{\alpha}\cap\Omega_{\beta})\times\mathbb D)$ and $\theta_{\alpha,\beta}(b)\in\mathbb R$. Moreover, it may be shown that $e^{i\theta_{\alpha,\beta}}$ is holomorphic and therefore locally constant.

Let $\tilde \pi:\tilde E\to B$ be a holomorphic fiber bundle with a fiber equal to $\tilde D=\{(z,w)\in\mathbb C^2:\ |z||w|<1\}$ and whose transition functions are defined in the following way  \begin{equation}\tilde \tau_{\alpha,\beta}(b,z,w)= \left(b,a_{\alpha,\beta}(b,zw)z, e^{iq\theta_{\alpha,\beta}(b)} \frac{w}{a_{\alpha,\beta}(b,zw)}\right).\end{equation} 

A direct computation allows us to observe that 
\begin{equation} E\in[(b,z,w)]\to[(b,z^p,w^q)]\in\tilde E \end{equation} is a well defined proper holomorphic mapping. Therefore, by a result of Narasimhan in \cite{Nar} (see also \cite{Gra}) the manifold $E$ is Stein if and only if $\tilde E$ is Stein. However, the Steinness of $\tilde E$ follows from the previous case.

$(c)$ Now we will show that $\mathcal D\times \mathcal C\in\mathfrak S$ for $\mathcal D\in\{\mathbb D,\mathbb D_*,\mathbb A(r)\}.$ 

Suppose that $E$ is a holomorphic fiber bundle with the fiber $\mathcal D\times\mathcal C$ and the Stein base $\Omega=\bigcup_{\alpha}\Omega_{\alpha}.$ By Theorem~\ref{wymaut}, any transition function of $E$ is of the form 
$$\tau_{\alpha,\beta}(b,\lambda,z)=(b,m_{\alpha,\beta}(b,\lambda),f_{\alpha,\beta}(b,\lambda,z)),\ (b,\lambda,z)\in(\Omega_{\alpha}\cap\Omega_{\beta})\times \mathcal D\times\mathcal C.$$ 

It follows from Theorem~\ref{wymaut} that $m_{\alpha,\beta}(x,\cdot)\in\aut(\mathcal D)$ for any $x\in\Omega_{\alpha}\cap\Omega_{\beta}$ (it may be shown that $m_{\alpha,\beta}(x,\cdot)$ does not depend on $x$ on connected components of $\Omega_{\alpha}\cap\Omega_{\beta}$).
 
Let $\tilde E$ be a holomorphic fiber bundle with the base $\Omega$ and the fiber $\mathcal D,$ whose transition functions $\tilde \tau_{\alpha,\beta}\in\aut((\Omega_{\alpha}\cap\Omega_{\beta})\times \mathcal D)$ are given by the formulas: 
\begin{equation}\tilde\tau_{\alpha,\beta}(b,\lambda)=(b,m_{\alpha,\beta}(b,\lambda)),\quad (b,\lambda)\in (\Omega_{\alpha}\cap\Omega_{\beta})\times \mathcal D.
\end{equation} By Theorem \ref{Mok0}, $\tilde E$ is Stein.

Now it is sufficient to observe that $$E\ni[(x,(\lambda,z))]\to[(x,\lambda)]\in\tilde E$$ forms a holomorphic fiber bundle with the base $\tilde E$ and the fiber equal to $\mathcal C.$ Using again Theorem \ref{Mok0} we get the Steinness of the bundle $E.$

$(d)$ To finish the proof of the theorem it suffices to show that $D=\{(z,w)\in\mathbb C_*\times\mathbb C:\ |z^p||w^q|<1\}\in\mathfrak S,$ where $p,q\geq 1$ are relatively prime natural numbers, $(p,q)\neq (1,1)$. We proceed similarly as in the case $(b).$ Namely, once again we aim at reducing situation to the already solved case $(c).$ 

Observe that $D=\{(z,w)\in\mathbb C_*\times\mathbb C:\ |z^p||w|<1\},$ $p\geq 2$ is algebraically equivalent to $\mathbb D\times \mathbb C_*$ (take the mapping $D\ni (z,w)\to (z^pw, w)\in \mathbb D\times \mathbb C_*$). Therefore we may assume that $q\geq 2.$

Suppose that $\pi:E\to B$ is a holomorphic fiber bundle. Then its transition functions must be of the form \begin{multline}\tau_{\alpha,\beta}(b,z,w)=(b,za^{1/p}_{\alpha,\beta}(b,z^pw^q),e^{i\theta_{\alpha,\beta}(b)}wa_{\alpha,\beta}^{-1/q}(b,z^pw^q)),\\ (b,z,w)\in(\Omega_{\alpha}\cap\Omega_{\beta})\times D,\end{multline} for some $a_{\alpha,\beta}\in\mathcal O^{*}((\Omega_{\alpha}\cap\Omega_{\beta})\times\mathbb D)$ and $\theta_{\alpha,\beta}(b)\in\mathbb R.$ It is seen that $e^{i\theta_{\alpha,\beta}}$ may be chosen to be constant on connected components of $\Omega_{\alpha}\cap\Omega_{\beta}.$ Therefore, putting $$\tilde \tau_{\alpha,\beta}(b,\lambda,z)=(b,e^{iq\theta_{\alpha,\beta}(b)}\lambda,za^{1/p}_{\alpha,\beta}(b,\lambda)),\quad (b,\lambda,z)\in(\Omega_{\alpha}\cap\Omega_{\beta})\times \mathbb D\times\mathbb C_*,$$ we obtain a holomorphic fiber bundle $\tilde E$ with the base $\Omega$ and the fiber $\mathbb D\times \mathbb C_*$ such that a mapping given by the formula
\begin{equation} E\ni[(x,z,w)]\rightarrow [(x,z^pw^q,z)]\in\tilde E \end{equation} is  proper and holomorphic.

The argument used in proof of $(b)$ finishes the proof of this case.
\end{proof}

\section{Irrational case}\label{sek2}
Geometry of Reinahrdt domains have been investigated in several papers (see \cite{Shi}, \cite{Shi2} and references contained there).

In \cite{ja} the author obtained a complete characterization of proper holomorphic mappings between Reinhardt domains of irrational type. This characterization is of the key importance for our considerations. Therefore, in Theorem~\ref{tkos} we collect the most  crucial part of it. 

Using Theorem~\ref{tkos} we shall obtain in the proof of Theorem~\ref{niewaut} a full description of the group of automorpisms of Reinahrdt domain of irrational type. As mentioned before, it turns out that this problem is connected with the famous Pell's equation.

Recall that any Reinhardt domain of irrational type is equivalent to a domain of the type (a), (b) or (c).
\begin{theorem}[\cite{ja}, see also \cite{Shi2}]\label{tkos}
(i) \begin{enumerate}\item[(a)] If $\alpha\in\mathbb R\setminus \mathbb Q,$ then the set of proper holomorphic mappings between $D_{\alpha,r}$ and $D_{\beta,R}$ is non-empty if and only if \begin{equation}\label{teza}\frac{\log R}{\log r}\in \mathbb Z + \beta \mathbb Z \quad \text{and} \quad \alpha\frac{\log R}{\log r}\in \mathbb Z + \beta\mathbb Z.\end{equation}

\item[(b)] Let $\alpha,\beta\in\mathbb R \setminus \mathbb Q,$ and $r,R>1$ be such that $\frac{\log R}{\log r}=k_1+l_1\beta$ and $\alpha\frac{\log R}{\log r}=k_2+l_2\beta$ for some integers $k_i,l_i,\ i=1,2.$ Then any proper holomorphic mapping $f:D_{\alpha,r}\rightarrow D_{\beta,R}$ is one of two following forms: \begin{equation}\label{glowne}\begin{cases}f(z)=(az_1^{k_1}z_2^{k_2},bz_1^{l_1}z_2^{l_2})\quad \text{or} \\ f(z)=(az_1^{-k_1}z_2^{-k_2},bz_1^{-l_1}z_2^{-l_2})\end{cases}\ z=(z_1,z_2)\in D_{\alpha,r},\end{equation} where $a,b\in\mathbb C$ satisfy the equality $|a||b|^{\beta}=1.$
\end{enumerate}

(ii) Let $\alpha,\beta\in\mathbb R\setminus \mathbb Q.$ The set of proper holomorphic mappings between $D_{\alpha}^*$ and $D_{\beta}^*$ is non-empty if and only if $\alpha=(k_2+\beta l_2)/(k_1+\beta l_1)$ for some $k_i,l_i\in\mathbb Z,\ i=1,2.$ Moreover, if $\alpha=(k_2+\beta l_2)/(k_1+\beta l_1),$ where $k_1+l_1\beta>0,$ then any proper holomorphic mapping $f:D_{\alpha}^* \rightarrow D_{\beta}^*$ is of the form \begin{equation}f(z_1,z_2)=(az_1^{k_1}z_2^{k_2},bz_1^{l_1}z_2^{l_2}),\quad (z_1,z_2)\in D_{\alpha}^*, \end{equation} where $a,b\in\mathbb C$ satisfy the equality $|a||b|^{\beta}=1.$

(iii) Let $\alpha,\beta\in\mathbb R \setminus \mathbb Q.$ 
\begin{enumerate}
\item[(a)] If $\alpha>0,\beta>0,$ then the set $\prop (D_{\alpha}, D_{\beta})$\footnote{$\prop(D,G)$ denotes the set of proper holomorphic mappings between $D$ and $G$.} is non-empty if and only if $\alpha=p\beta$ for some $p\in\mathbb Q_{>0}.$ In this case all proper holomorphic mappings between $D_{\alpha}$ and $D_{\beta}$ are of the form \begin{equation}(z_1,z_2)\to(z_1^k,z_2^l),\end{equation} where $k,l$ are positive integers such that $p=\frac{l}{k}.$

\item[(b)] If $\alpha<0,\beta<0,$ then the set $\prop(D_{\alpha},D_{\beta})$ is non-empty if and only if $\alpha=p_1+p_2\beta$ for some rational $p_1,p_2,$ $p_2\neq0.$ In this case proper holomorphic mappings between $D_{\alpha}$ and $D_{\beta}$ are of the form \begin{equation} (z_1,z_2) \to  (z_1^{k_1}z_2^{k_2},z_2^l) ,\end{equation} where $k_1,k_2,l,$ $k_1>0,$ are integers such that $p_1=\frac{k_2}{k_1},$ i  $p_2=\frac{l}{k_1}.$

\item[(c)] If $\alpha\beta<0,$ then there is no proper holomorphic mapping between $D_{\alpha}$ and $D_{\beta}.$
\end{enumerate}
\end{theorem}

Our aim in this part of the paper is to show the following
\begin{theorem}\label{niewaut} Let $\alpha\in\mathbb R\setminus\mathbb Q.$ Then \begin{enumerate}

\item[(a)] $D_{\alpha}^*\notin\mathfrak S $ if and only if $\alpha=p\pm\sqrt{q}$ for some $p,q\in\mathbb Q,\ q>0,$
\item[(b)] $D_{\alpha}\in\mathfrak S $ for any $\alpha\in\mathbb R\setminus\mathbb Q,$
\item[(c)] $D_{\alpha,r}\in\mathfrak S $ for any $\alpha\in\mathbb R\setminus\mathbb Q,\ r>1.$
\end{enumerate}
\end{theorem}

\begin{proof}
(a) It follows from Theorem~\ref{tkos} that every automorphism of the domain $D^*_{\alpha}$ is of the form 
\begin{equation}\label{aut}(z_1,z_2)\to (az_1^{k_1}z_2^{k_2}, bz_1^{l_1}z_2^{l_2}),
\end{equation} where $|a||b|^{\alpha}=1$ and $k_1,k_2,l_1,l_2\in\mathbb Z$ are such that $k_1l_2-k_2l_1=\pm1,$ $\alpha=(k_2+l_2\alpha)/(k_1+l_1\alpha)$ and $k_1+l_1\alpha>0.$

These conditions imply that if $\alpha$ is not of the form $p\pm\sqrt{q}$ for any rational $p,q,$ then $k_1=l_2=1$ and $k_2=l_1=0.$ In other words any automorphism of $D_{\alpha}^*$ is of the form $(z_1,z_2)\to (a z_1,b z_2),\ a,b\in\mathbb C,\ |a||b|^{\alpha}=1.$ Therefore, the plurisubharmonic function $u$ given by the formula $$u(z_1,z_2)=\max\left\{\frac{1}{1-|z_1||z_2|^{\alpha}},\log|z_1|,-\log|z_1|,\log|z_2|,-\log|z_2|\right\},$$
 $(z_1,z_2)\in D_{\alpha}^*,$ satisfies the Stehl\'e criterion, hence $D_{\alpha}^*\in\mathfrak S .$

Suppose now that $\alpha=\frac{p}{n}\pm\sqrt{\frac{q}{n^2}}$ for some $p,q\in\mathbb Z,\ q>0,\ n\in\mathbb N.$ 

We are looking for $k_i,l_i\in\mathbb Z,\ k_1+\alpha l_1\geq0$ such that ($\ref{aut}$) defines an automorphism of $D_{\alpha}^*.$ Then the following equations are satisfied
\begin{equation}\label{wyz}\begin{cases}k_2=l_1(\frac{q}{n^2}-\frac{p^2}{n^2}),\\ l_2=k_1+2\frac{p}{n}l_1,\end{cases}
\quad \text{and}\quad |k_1l_2-k_2l_1|=1.\end{equation} 

Consider the so-called \textit{Pell's equation} of the following form:
\begin{equation}\label{pell} n^2q=\frac{x^2-1}{y^2}.
\end{equation} It was shown by Lagrange that (\ref{pell}) has infinitely many integer solutions (recall that $q$ is not a square of a natural number). Let $x,y,\ x,y>0,$ denote the arbitrary natural solution of this equation. Put \begin{equation}l_1=n^2y,\quad k_1=x-pny,\quad k_2=y(q-p^2),\quad \text{and}\quad l_2=x+pny.\end{equation}  

It is a direct consequence of (\ref{pell}) that $x\pm ny\sqrt{q}>0,$ hence $k_1+\alpha l_1>0.$ An easy computation shows that each condition in $(\ref{wyz})$ is satisfied by such chosen integers $k_i,l_i,\ i=1,2.$

Moreover $\tr\left(%
\begin{array}{cc}
  k_1 & k_2 \\
  l_1 & l_2 \\
\end{array}%
\right)>2$ so it follows from Theorem~1 in \cite{Zaf} that $D\not\in\mathfrak S.$

(b) It follows from Theorem~\ref{tkos} that any automorphism of $D_{\alpha}$ is elementary algebraic. The automorphisms must also preserve the axis $\{0\}\times\mathbb C_*$ (when $\alpha>0$ additionally the axis $\mathbb C\times\{0\}$ is also preserved). From this piece of information one may conclude that any automorphism of domain $D_{\alpha}$ is of the form \begin{equation}(z_1,z_2)\to(az_1,bz_2),\ (z_1,z_2)\in D_{\alpha},\end{equation}
where $a,b\in\mathbb C,\ |a||b|^{\alpha}=1.$ Therefore the functions $u_+,u_-$ given by the formulas \begin{align*}u_+(z_1,z_2)&=\max\left\{\frac{1}{1-|z_1||z_2|^{\alpha}},\log|z_1|,\log|z_2|\right\},\quad \text{when}\ \alpha>0,\\
\text{and}\ u_-(z_1,z_2)&=\max\left\{u_+(z_1,z_2),-\log|z_2|\right\},\quad \text{when}\ \alpha<0,\end{align*} satisfy the criterion of Stehl\'{e}.

(c) Similarly as before we easily find from Theorem~\ref{tkos} that the group of automorphisms of $D_{\alpha,r}$ consists of the mappings of the form: 
\begin{align}(z_1,z_2)\to(az_1^{\epsilon},bz_2^{\epsilon}),\ (z_1,z_2)\in D_{\alpha,r},\end{align} where $\epsilon=\pm1$ and $a,b\in\mathbb C,\ |a||b|^{\alpha}=1.$ Thus the function $$u(z_1,z_2)=\max\left\{\frac{|z_1||z_2|^{\alpha}}{r|z_1||z_2|^{\alpha}-1},\frac{1}{r-|z_1||z_2|^{\alpha}},\log^2|z_1|,\log^2|z_2|\right\},$$ $(z_1,z_2)\in D_{\alpha},$ satisfies the criterion of Stehl\'{e}.
\end{proof}

\section{The case when $D\cap \mathbb C_*^2$ is hyperbolic}\label{sek3}
For a pseudoconvex Reinhardt domain $D$ in $\mathbb C^n$ let $I(D)$ denote the set of $i=1,\ldots,n,$ for which the intersection $(\mathbb C^{i-1}\times\{0\} \times\mathbb C^{n-i})\cap D$ is not hyperbolic (viewed as a domain in $\mathbb C^{n-1}$). Put
\begin{equation}D^{hyp}=D\setminus \big(\bigcup_{i\in I(D)} (\mathbb C^{i-1}\times\{0\} \times\mathbb C^{n-i})\big).
\end{equation}

With a Reinhardt domain $D$ we associate the following constant:
$$t(D):=\ \textit{the number of}\ i=1,\ldots,n,\ \textit{such that}\ D\cap (\mathbb C^{i-1}\times\{0\} \times\mathbb C^{n-i})\neq\emptyset.$$

\begin{example}
 Put $D_1:=\{(z,w)\in\mathbb C\times\mathbb C:\ |z|<1,\ |zw|<1\}$ and $D_1:=\{(z,w)\in\mathbb C_*\times\mathbb C:\ |z|<1,\ |zw|<1\}.$ Then both $D_1$ and $D_2$ are non-hyperbolic. Observe that $t(D_1^{hyp})=1$ and $t(D_2^{hyp})=0.$
\end{example}

For a non-hyperbolic pseudoconvex Reinhardt domain $D\subset\mathbb C^2$ such that $D^{hyp}$ is hyperbolic the number $t(D^{hyp})$ is equal either to $1$ or to $0.$ 

First we shall focus our attention on the case $t(D^{hyp})=1.$ We start with the following
\begin{lemma}\label{lem} Let $D$ be a non-hyperbolic, pseudoconvex Reinhardt domain in $\mathbb C^2.$ Assume that the logarithmic image of $D$ contains no affine line and the group of automorphisms of $D^{hyp}$ is compact. Then $\aut(D)$ is also compact.
\end{lemma}
\begin{proof} Take any sequence $(\varphi_n)_n\subset\aut(D).$ Since for any $\varphi\in\aut(D)$ the restriction $\varphi|_{D^{hyp}}$ is an automorphism of $D^{hyp}$ (see e.g. \cite{ja2}, Theorem 8), we may assume that $({\varphi_n}|_{_{D^{hyp}}})_n$ is convergent locally uniformly on $D^{hyp}.$ Applying the Cauchy's formula we infer that $(\varphi_n)_n$ is convergent to some holomorphic function on $D.$ Repeating the above argument for the sequence $(\varphi_n^{-1})_n$ immediately gives the desired result.
\end{proof}

The problem of a characterization of automorphism groups of bounded Reinhardt domains was studied in \cite{Shi} (see also \cite{isa-kra}, \cite{Isa} and papers quoted there).

The results obtained in \cite{Shi} and \cite{Kru} together with remarks from \cite{pzwo} lead us to the description of pseudoconvex hyperbolic Reinhardt domains with $t=1$ and non-compact automorphism group. For our future use we recall here the version formulated in \cite{pzwo}.
\begin{theorem}[see \cite{Kru},\cite{Shi} and \cite{pzwo}, Theorem 4]\label{pz} Let $D$ be a hyperbolic, pseudoconvex Reinhardt domain with $t(D)=1.$ Then $\aut(D)$ is non-compact if and only if $D$ is algebraically equivalent to one of the following domains:
\begin{enumerate}[(a)]
\item $\mathbb D\times \mathbb A(r,1),$ where $0\leq r<1.$ In this case the group of automorphisms consists of the mappings of the form $D\ni(z_1,z_2) \to(a(z_1),b(z_2))\in D,$ where $a\in\aut(\mathbb D)$ and $b\in\aut(\mathbb A(r,1)).$
\item $\{(z_1,z_2)\in\mathbb C^2:\ |z_1|<1,\ 0<|z_2|<(1-|z_1|^2)^{p/2}\},\quad p>0.$ In this case the group of automorphisms consists of the mappings of the form $D\ni(z_1,z_2) \to(\alpha\frac{z_1-\beta}{1-\overline{\beta}z_1},\gamma\frac{(1-|\beta|^2)^{\frac{p}{2}}}{(1-\overline{\beta}z_1)^p}z_2)\in D,$ where $|\alpha|=|\gamma|=1,\ \beta\in\mathbb C.$
\item $\{(z_1,z_2)\in\mathbb C^2:\ 0<|z_2|<\exp(-|z_1|^2)\}.$ In this case the group of automorphisms consists of the mappings of the form $D\ni(z_1,z_2) \to(\alpha z_1+\beta,\gamma\exp(-2\alpha\overline{\beta}z_1-|\beta|^2)z_2)\in D,$ where $|\alpha|=|\gamma|=1,\ \beta\in\mathbb C.$
\end{enumerate}

The following lemma will be also needed:
\begin{lemma}[see \cite{ja2}, Theorem~2]\label{restr} Let $D,$ $G$ be pseudoconvex Reinhardt domains in $\mathbb C^2$. Let $\varphi: D\to G$ be a proper holomorphic mapping. Assume that the logarithmic images of domains $D$ and $G$ contain no affine lines. Then $$\varphi(D^{hyp})\subset G^{hyp}$$ and the restriction $$\varphi|_{D^{hyp}}: D^{hyp} \to G^{hyp}$$ is proper.
\end{lemma}

\end{theorem}

\begin{theorem}\label{14} Let $D\subset\mathbb C^2$ be a pseudoconvex non-hyperbolic Reinhardt domain with $t(D^{hyp})=1.$ Assume additionally that the logarithmic image of the domain $D^{hyp}$ contains no affine lines. 
Then the group of automorphism of domain $D$ is compact. In particular, the answer to the Serre problem for the domain $D$ is positive.
\end{theorem}
\begin{proof} In view of Lemma~\ref{lem} it suffices to prove the theorem for domain $D$ under the additional assumption that $\aut(D^{hyp})$ is non-compact. Then it is a consequence of Theorem~\ref{pz} that up to a multiplying and a permutation of components $D$ must have the form:
\begin{enumerate}[(a)]
\item $\{(z_1,z_2)\in\mathbb C^2:\ |z_1z_2^k|<1,\ |z_2|<1\},$ where $k\in\mathbb Z_{>0},$
\item $\{(z_1,z_2)\in\mathbb C^2:\ |z_1z_2^k|<1,\ |z_2|<(1-|z_1z_2^k|^2)^{p/2}\}$, for some $p>0,\ k\in\mathbb Z_{>0},$
\item $\{(z_1,z_2)\in\mathbb C^2:\ |z_2|<\exp(-|z_1z_2^k|^2))\},\ k\in\mathbb Z_{>0}.$
\end{enumerate} 

Moreover (see Lemma~\ref{restr}) any automorphism of the domain of one of the forms presented above preserves the axis $\mathbb C\times\{0\}.$ This fact and Theorem~\ref{pz} lead to the statement that the group of automorphisms of the domain D in all cases (a),(b)  and (c) consists of the mappings of the form $$D\ni(z_1,z_2)\to(az_1,bz_2)\in D,$$ where $|a|=|b|=1.$ Therefore $\aut(D)$ is compact in all cases; a contradiction.
\end{proof}

Let us pass to the remaining case $t(D^{hyp})=0$. An important role in our approach is played by the following result
\begin{theorem}[see \cite{Kru}, \cite{Shi}]\label{sh}$\aut(D)=\aut_{alg}(D)$ for any pseudoconvex hyperbolic Reinhardt domain $D\subset\mathbb C_*^n.$
\end{theorem}

\begin{theorem}\label{16} Let $D\subset\mathbb C^2$ be a non-hyperbolic pseudoconvex Reinhardt domain such that $D^{hyp}$ is hyperbolic and $t(D^{hyp})=0.$ Then $D\in\mathfrak S.$

Moreover, the group of automorphism of $D$ is non-compact if and only if, up to a permutation on components, $D$ is contained in $\mathbb C\times\mathbb C_*,$ $\{0\}\times\mathbb C_*\subset D$ and \begin{equation}\label{niec}\log D=\{(t,s):\ t<\psi(s),\ s\in\mathbb R\},\end{equation} where $\psi:\mathbb R\to \mathbb R$ is a concave function satisfying the property $\psi(\beta+ s)-\psi(s)=\alpha+ks,\ s\in\mathbb R$ for some $ \alpha,\beta\in\mathbb R,\ \beta\neq 0,\ k\in\mathbb Z_*.$
\end{theorem}
\begin{proof}Using the inclusion $\aut(D)_{|{D^{hyp}}}\subset\aut(D^{hyp})$ (Lemma~\ref{restr}) and Theorem~\ref{sh} we state that any automorphism of $D$ must be algebraic. Moreover, at least one of the axis $\{0\}\times\mathbb C_*,\ \mathbb C_*\times\{0\}$ must be contained in D (otherwise $D$ would be hyperbolic). 

Suppose that both axes are contained in $D.$ Then, since any automorphism of $D$ maps the axes onto the axes (use Lemma~\ref{restr}), we see that the group $\aut(D)$ consists of the mappings of the form
\begin{align}\label{oa}D\ni(z_1,z_2)\rightarrow (az_1,bz_2)\in D\quad \text{or}\\\label{oa2} \quad D\ni(z_1,z_2)\rightarrow (az_2,bz_1)\in D.
\end{align} for suitable $a,b\in\mathbb C_*.$ We shall show that $|a|=|b|=1$ in the case (\ref{oa}). Moreover, we shall prove that there is $R>0$ such that for any automorphism satisfying (\ref{oa2}) $|a|=R$ and $|b|=1/R.$ This in particular means that $\aut(D)$ is compact.

Let us take any $\varphi\in\aut(D)$ of the form~(\ref{oa}).

Assume a contrary, i.e. $(\log|a|,\log|b|) \neq(0,0).$ For $n\in\mathbb N$ put $\varphi^{(n)}=\varphi\circ\ldots\circ\varphi\in\aut(D)$ and $ \varphi^{(-n)}=(\varphi^{-1})^{(n)}.$ Since $\varphi^{(n)}(z_1,z_2)=(a^nz_1,b^nz_2)\in D$ for any $(z_1,z_2)\in D$ passing to the logarithmic image of the domain $D$ easily shows that $\mathbb R(\log|a|,\log|b|)+\log D\subset\log D.$ This is in a contradiction with the hyperbolicity of $D^{hyp}.$

Assume that $\varphi:(z_1,z_2)\rightarrow(az_2,bz_1)$ is an automorphism of the domain $D.$ Repeating the previous reasoning applied to an automorphism $\varphi^{(2)}(z_1,z_2)=(abz_1,abz_2)$ we immediately see that $|ab|=1.$ Therefore it suffices to show that for any other automorphism $\tilde\varphi:(z_1,z_2)\rightarrow(a_1z_2,b_1z_1)$ of $D$ the following relations holds: $|a_1|=|a|$ and $|b_1|=|b|.$ However, these relations also follows from the above reasoning applied to $\tilde\varphi\circ\varphi(z_1,z_2)=(a_1bz_1,ab_1z_2).$ 

This finishes the proof of the compactness of $\aut(D)$ in this case.

Suppose that only one axis is contained in $D$, e.g. $\{0\}\times\mathbb C_*\subset D.$ We use the idea applied by P.~Pflug and W.~Zwonek in~\cite{pzwo}. 

Assume that the group $\aut(D)$ is non-compact.
Note that $D\cap(\mathbb C\times\{0\})=\emptyset.$ As before, using the fact that automorphisms preserve the axis (Lemma~\ref{restr}), we see that $\aut(D)$ consists of the mappings of the form 
\begin{equation}\label{a8}\Phi=\Phi_{a,b,k,\epsilon}:(z_1,z_2)\to (az_1z_2^k,bz_2^{\epsilon})
\end{equation} for some $a,b\in\mathbb C_*,$ $\epsilon=\pm1$ and $k\in\mathbb Z.$ 

(\dag) We shall show that there exists an automorphism $\Phi_{a,b,k,\epsilon}$ of the domain $D$ with $|b|\neq1$ and $\epsilon=1$ (then also $k\neq 0$). 

First observe that there is an automorphism of the domain $D$ of the form (\ref{a8}) with $k\neq0$. To show it repeat the argument from the first part of the proof to observe that $|a|=|b|=1$ for any automorphism $\Phi_{a,b,0,1}.$ Moreover, there exists $R>0$ such that $|a|=R,$ $|b|=1/R$ for any automorphism $\Phi_{a,b,0,-1}.$ So if there were no automorphisms of the form (\ref{a8}) with $k\neq 0$, then the group $\aut(D)$ would be compact.

If there is an automorphism $\Phi$ of $D$ of the form (\ref{a8}) with $\epsilon=1$ and $k\neq 0,$ we are done. Actually, $\Phi^{(n)}(z_1,z_2)=(a^nb^{kn(n-1)/2}z_1z_2^{nk},b^n z_2)$ and using the hyperbolicity of $D^{hyp}$ once again, we find that $|b|\neq 1.$ 

Thus, to prove (\dag) it suffices to consider the case when any automorphism of the form (\ref{a8}) with $k\neq 0$ satisfies $\epsilon=-1.$ Since $\Phi_{a,b,k,-1}^{(2n)}(z_1,z_2)=((a^2b^k)^nz_1,z_2),$ $n\in\mathbb Z,$ we see that $|a^2b^k|=1.$ So there is an automorphism $\Phi_{\tilde a, \tilde b, \tilde k,-1}$ such that $|\tilde b|\neq |b|$ or $\tilde k\neq k$ (as $\aut(D)$ is non-compact). Let us compute $$\Phi_{a,b,k,-1}\circ\Phi_{\tilde a,\tilde b,\tilde k,-1}=\Phi_{a\tilde a\tilde b^k,\frac{b}{\tilde b},\tilde k-k,1}.$$ If $|\tilde b|\neq |b|,$ then the statement is clear. Otherwise $k\neq \tilde k$. But this is in a contradiction with the assumption that any automorphism of the form (\ref{a8}) with $k\neq 0$ satisfies $\epsilon=-1.$

So we have shown (\dag).

Properties of pseudoconvex Reinhardt domains imply that for any $s\in\mathbb R$ there is (exactly one) $\psi(s)\in\mathbb R$ such that $(\psi(s),s)\in\partial \log D.$ Moreover, it is an immediate consequence of  the inclusion $\{0\}\times\mathbb C_*\subset D$ and convexity of $\log D$ that the function $\psi(\cdot)$ is concave. Put $$v(z_1,z_2)=\log|z_1|-\psi(\log|z_2|).$$

For fixed $(z_1,z_2)\in D\cap \mathbb C_*^2$ put $(t,s)=(\log|z_1|,\log|z_2|).$ Let $\Phi=\Phi_{a',b',k',\epsilon}$ be an arbitrary automorphism given by the formula (\ref{a8}). 

Denote $(t',s')=(\log|\Phi_1(z)|,\log|\Phi_2(z)|).$ It is seen that $(t',s')=(\log|a'|+t+k's,\log|b'|+\epsilon s).$ 

Since $\Phi((\partial D)\cap\mathbb C_*^2))=(\partial \Phi(D))\cap\mathbb C_*^2,$ we see that \begin{equation}\label{t(s)}(\psi(s'),s')=(\log|a'|+\psi(s)+k's,\log|b'|+\epsilon s).\end{equation} Therefore $\psi(s')-t'=\psi(s)-t.$ In other words
\begin{equation}\label{v}
 v\circ\Phi=v\quad\text{for any}\quad \Phi\in\aut(D).
\end{equation}

Define $$u(z_1,z_2)=\max\left\{\log|z_2|,-\log|z_2|,-(v(z_1,z_2))^{-1}\right\},\quad (z_1,z_2)\in D.$$ A direct calculation (approximate and compute the Levi form) shows that $u$ is plurisubharmonic. Moreover, $u$ is an exhausting function for $D.$ In view of (\ref{v}) $u$ satisfies the criterion of Stehl\'e.

Using (\dag) and the property (\ref{t(s)}) once again we find that $\psi(\log |b|+ s)-\psi(s) =\log|a|+ks.$ From this we obtain a desired properties of the function $\psi.$

Conversely, having a domain $D\subset\mathbb C\times\mathbb C_*$ satisfying (\ref{niec}) one may easily see that the mapping $\Phi$ given by the formula~(\ref{a8}) with $b=e^{\beta},\ a=e^{\alpha},\ \epsilon=1$ is an automorphism of the domain $D.$ The investigation of $\Phi^{(n)}$ immediately proves a non-compactness of $\aut(D).$
\end{proof}

\section{$\mathbb C\times\mathbb C_*$}
In 1977 a negative answer to the Serre problem was given by H. Skoda who proved that $\mathbb C^2\notin\mathfrak S.$ This construction was improved in \cite{Dem76} by J.P. Demailly who proved that polynomial automorphisms of $\mathbb C^2$ may serve as the transition function.  
Later in~\cite{Dem} J.P. Demailly constructed a counterexample to the Serre problem with a plane or a disc as a base. 

Let us recall here this construction. The base $\Omega$ is a domain containing $3\mathbb D.$ Put $\Omega_0=\Omega\setminus\{-1,1\},\ \Omega_1=\Omega_0\cup\{1\},\ \Omega_2=\Omega_0\cup\{-1\}.$ The transition functions $\tau_{i,j}:\Omega_i\times\mathbb C^2\to\Omega_j\times\mathbb C^2,\ i\neq j$ of the fiber bundle $E$ are defined as follows:
\begin{align*}&\tau_{0,1}(x;z_1,z_2)=(x;z_1,z_2 \exp(z_1u(x)),\\
&\tau_{0,2}(x;z_1,z_2)=(x,z_1\exp(z_2u(x)),z_2)\end{align*} and $\tau_{1,2}=\tau_{0,1}^{-1}\circ\tau_{0,2},$ where $u(x)=\exp(\frac{1}{x^2-1}).$ It is clear that any plurisubharmonic function $V$ on $E$ induces plurisubharmonic $V_j$ such that $V_j=V_k\circ\tau_{k,j}.$ The idea of the proof relied upon comparing the maximum of the functions $V_j$ over the polidiscs $\frac{1}{2}\mathbb D\times (r\mathbb D)^2,\ r>>1.$ More precisely, it was shown that \begin{equation}\label{m}M(V_0,\frac{1}{2}\mathbb D,\exp(r/32))\leq M(V_0,\frac{1}{2}\mathbb D,\exp(\log^3r))+C,\quad r>>1,\end{equation} where $M(V,\omega,r)=\max_{\omega\times (r\mathbb D)^2}V$ and the constant $C$ does not depend on $r.$ The key role was played by the logarithmic convexity of the functions $$(\rho,r)\to M(V,\rho\mathbb D,r),\quad V\in\mathcal{PSH}(\Omega\times\mathbb C^2).$$ 

Direct computations allows us to obtain logarithmic convexity of the function $(\rho,r)\to\max_{\rho\mathbb D\times r\mathbb D\times\mathbb A(r)}\tilde V$ for any $\tilde V\in\mathcal{PSH}(\Omega\times\mathbb C\times\mathbb C_*).$
Therefore considering $\tilde M$ instead of $M,$ where $\tilde M(V,\omega,r)=\max_{\omega\times r\mathbb D\times\mathbb A(r)}V$ and repeating the reasoning from Demailly's paper we may replace $M$ by $\tilde M$ in the inequality (\ref{m}). This, together with mentioned above logarithmic convexity of $\tilde M$ immediately shows that $\mathbb C\times\mathbb C_*$ does not belong to $\mathfrak S.$

\section{Proofs of the main theorems}\label{pomt}
\begin{proof}[Proof of Theorem~\ref{main}]
The result follows from Theorem~\ref{dc}, Theorem~\ref{niewaut}, Theorem~\ref{14} and Theorem~\ref{16}.
\end{proof}

\begin{proof}[Proof of Theorem~\ref{compact}]
It is clear that if $\log D$ contains an affine line, then the group $\aut(D)$ is non-compact.

In the case when $\log D$ contains no affine lines (i.e. $D^{hyp}$ hyperbolic) the result follows immediately from Theorems~\ref{14} and~\ref{16}.
\end{proof}


\begin{thebibliography}{999999}
\bibitem[Coe-Loeb]{Coe}\textsc{G.~Coeur\'e and J.-J.~Loeb}  \textit{A counterexample to the Serre problem with a bounded domain in $\mathbb C^2$ as fiber}, Ann. Math., 122 (1985), 329-334.
\bibitem[Dem1]{Dem76}\textsc{J.P.~Demailly,}  \textit{Diff\'erents exemples de fibr\'e holomorphes non de Stein}, S\'eminaire P. Lelong - H. Skoda, 1976-77, 15-41, Lecture Notes in Math., 694.
\bibitem[Dem2]{Dem}\textsc{J.P.~Demailly,}  \textit{Un exemple de fibr\'e holomorphe non de Stein
\'a fibre $C^2$ au-dessus du disque ou du plan}, S\'eminaire P. Lelong, P. Dolbeault, H. Skoda (Analyse) \textbf{24} (1983/84)
Lecture Notes in Math. 1198, Springer, 88-97.
\bibitem[Edi-Zwo]{prop} \textsc{A.~Edigarian and W.~Zwonek,} \textit{Proper holomorphic mappings in some class of unbounded domains}, Kodai Math J. \textbf{22} (1999), 305--312.
\bibitem[Gra]{Gra}\textsc{H.~Grauert} \textit{Charakterisierung der holomorph-vollst\"{a}ndigen komplexen Räume}, Math. Ann. 129 (1955), 233-259.
\bibitem[H\"{o}r]{Her}\textsc{L.~H\"ormander} \textit{An introduction to complex analysis in several variables}, third edition, North-Holland Mathematical Library, 7. North-Holland Publishing Co., Amsterdam, 1990.
\bibitem[Isa-Kra]{isa-kra} \textsc{A.V.~Isaev and S.G.~Krantz}, \textit{Domains with non-compact automorphism group: a survey,} Adv. Math. \textbf{146} (1999), 1–38.
\bibitem[Isa-Kru]{Isa} \textsc{A.V.~Isaev and N.G.~Kruzhilin,} \textit{Proper holomorphic Maps between Reinhardt domains in
$\mathbb C^2$}, Michigan Math. \textbf{54} (2006), 33--64.
\bibitem[Jar-Pfl]{Jar-Pfl} \textsc{M.~Jarnicki and P.~Pflug,} \textit{First Steps in Several Complex Variables: Reinhardt domains,} EMS Textbooks in Mathematics. European Mathematical Society (EMS), Z\"urich, 2008.
\bibitem[Kos1]{ja}\textsc{\L.~Kosi\'nski,}  \textit{Proper holomorphic mappings in the special class of Reinhardt domains},
Ann. Polon. Math. \textbf{92} (2007), 285-297.
\bibitem[Kos2]{ja2}\textsc{\L.~Kosi\'nski,}  \textit{Proper holomorphic mapppings between Reinhardt domains in $\mathbb C^2$},
Mich. Math. Journal. \textbf{58} (2009), no. 3, 711--721.
\bibitem[Kru]{Kru} \textsc{N.G.~Kruzhilin}, \textit{Holomorphic equivalence of hyperbolic Reinhardt domains}, Math. USSR Izv., \textbf{32} (1989), 15-38.
\bibitem[Mok1]{Mok0} \textsc{N.~Mok}, \textit{Le probl\`{e}me de Serre pour les surfaces de Riemann}, C. R. Acad. Sci. Paris S\'er. A-B 290 (1980) \textbf{4}, A179-A180.
\bibitem[Mok2]{Mok} \textsc{N.~Mok}, \textit{The Serre problem on Riemann surfaces}, Math. Ann., (258) (1981), 145-168.
\bibitem[Nar]{Nar} \textsc{R.~Narasimhan}, \textit{A note on Stein spaces and their normalizations}, Ann. Sc. Norm. Sup. Pisa, 16 (1962), 327-333.
\bibitem[Oel-Zaf]{Oel} \textsc{K.~Oeljeklaus and D.~Zaffran}, \textit{Steinness of bundles with fiber a Reinhardt bounded domain}, Bull. Soc. math. France \textbf{134} 4 (2006), 451-473.
\bibitem[Pfl-Zwo]{pzwo} \textsc{P.~Pflug and W.~Zwonek,} \textit{The Serre problem with Reinhardt fibers},
Annales de l'Institut Fourier  \textbf{54} (2004), 129-146.
\bibitem[Shi1]{Shi} \textsc{S.~Shimizu}, \textit{Automorphisms and equivalence of bounded Reinhardt domains not containing the origin}, Tohoku Math. J., (40) \textbf{1} (1988), 119-152.
\bibitem[Shi2]{Shi2} \textsc{S.~Shimizu}, \textit{Holomorphic equivalence
problem for a certain class of unbounded Reinhardt domains in $\mathbb C^2$}, II, Kodai Math. J, \textbf{15} (1992), 430-444.
\bibitem[Sko]{Sko} \textsc{H.~Skoda}, \textit{Fibr\'e holomorphes \'a base fibre et \'a fibre de Stein,} Invent. Math. \textbf{43} (1977), 97-107.
\bibitem[Sol]{Sol} \textsc{P.A.~Soldatkin}, \textit{Holomorphic equivalence of Reinhardt domains in $\mathbb C^2$}. (Russian) Izv. Ross. Akad. Nauk Ser. Mat. \textbf{66} (2002), no. 6, 187--222; translation in Izv. Math. \textbf{66} (2002), no. 6, 1271--1304.
\bibitem[Ste]{Stel} \textsc{J.-L.~Stehl\'e}, \textit{Fonctions plurisousharmoniques et convexit\'e holomorphe de certaines fibr\'es analytiques}, In S\'eminaire Pierre Lelong (Analyse), Ann\'ee 1973/74 (ed. P. Lelong, P. Doulbeaut, H. Skoda), Lecture Notes in Mathematics 474. Springer, Berlin (1975)m 155-179.
\bibitem[Zaf]{Zaf} \textsc{D.~Zaffran}, \textit{Holomorphic functions on bundles over annuli}, Math. Ann. 341 (2008), no. \textbf{4}, 717-733.
\bibitem[Zwo]{Zwo} \textsc{W.~Zwonek,} \textit{On hyperbolicity of pseudoconvex Reinhardt domains},
Archiv der Mathematik \textbf{72} (1999), 304-314.
\end{thebibliography}
\end{document}